\documentclass[11pt]{article}
\usepackage[margin=1in]{geometry}
\usepackage{hyperref}
\usepackage{graphicx,amsmath,color}
\usepackage{amssymb,amsthm,mathrsfs}
\usepackage{algorithmic,algorithm}
\usepackage{mathtools}
\usepackage{subfig}
\usepackage{todonotes}
\usepackage{cleveref}
\DeclareGraphicsExtensions{.pdf}

\graphicspath{{Fig/}}

\newcommand*\Laplace{\mathop{}\!\mathbin\bigtriangleup}

\newcommand{\bq}{\begin{equation}}
\newcommand{\eq}{\end{equation}}

\definecolor{rev1}{rgb}{0.79,0.15,0.06}

\newcommand{\G}{\mathcal{G}}
\newcommand{\bO}{\mathcal{O}}

\newcommand{\cL}{\mathcal{L}}
\newcommand{\cM}{\mathcal{M}}

\DeclareMathOperator*{\argmin}{argmin}

\theoremstyle{remark}

\title{Anderson Acceleration for Seismic Inversion}
\author{Yunan Yang\footnote{Courant Institute of Mathematical Sciences, New York University, New York, NY 10012. (yunan.yang@nyu.edu)}
}

\begin{document}

\maketitle

\begin{abstract}
The state-of-art seismic imaging techniques treat inversion tasks such as FWI and LSRTM as PDE-constrained optimization problems. Due to the large-scale nature, gradient-based optimization algorithms are preferred in practice to update the model iteratively. Higher-order methods converge in fewer iterations but often require higher computational costs, more line search steps, and bigger memory storage. A balance among these aspects has to be considered. We propose using Anderson acceleration (AA), a popular strategy to speed up the convergence of fixed-point iterations, to accelerate the steepest descent algorithm, which we innovatively treat as a fixed-point iteration. Independent of the dimensionality of the unknown parameters, the computational cost of implementing the method can be reduced to an extremely low-dimensional least-squares problem. The cost can be further reduced by a low-rank update. We discuss the theoretical connections and the differences between AA and other well-known optimization methods such as L-BFGS and the restarted GMRES and compare their computational cost and memory demand. Numerical examples of FWI and LSRTM applied to the Marmousi benchmark demonstrate the acceleration effects of AA. Compared with the steepest descent method, AA can achieve fast convergence and provide competitive results with some quasi-Newton methods, making it an attractive optimization strategy for seismic inversion. 
\end{abstract}

\section{Introduction}
The fast growth of computational power popularizes numerous techniques that utilize the full wavefields in seismic imaging~\cite{tarantola1982}.  In particular, the full-waveform inversion (FWI)~\cite{virieux2009overview} and the least-squares reverse-time-migration (LSRTM)~\cite{Dai2013} aim to reconstruct the subsurface properties such as the wave velocity and the material density by minimizing an objective function that measures the discrepancy between the synthetic data and the observed data. Iterative optimization algorithms are then applied to find the optimal solution~\cite{metivier2012truncated}. 

For local optimization, the descent direction depends on the gradient and the Hessian information of the objective function with respect to the model parameters. Theoretically, the step size along the descent direction should be determined by line search to guarantee a sufficient decrease in the objective function and avoid overshooting. However, the process of backtracking line search could incur a considerable amount of extra wave modeling. Sometimes, to reduce the computational cost of line search and avoid overshooting, tiny fixed step size is preferred instead, but it slows down the convergence. Similarly, Newton's method is not widely used in practical seismic inversions due to the cost of calculating and storing the Hessian matrix, despite known to offer a quadratic convergence rate. For large-scale optimization problems such as seismic inversion, a better rate of convergence often comes at the cost of memory and computing power. The best practice is continuously a balance between the two aspects.

In the last two decades, the so-called Anderson acceleration (AA) has been widely used in several applied fields for problems that can be solved by a fixed-point iteration. The application of AA includes flow problems~\cite{pollock2018anderson}, solving nonlinear radiation-diffusion equations~\cite{an2017anderson} and~wave propagation~\cite{yang2020anderson}.
It is closely related to Pulay mixing~\cite{pulay1980convergence} and DIIS (direct inversion on the iterative subspace)~\cite{kudin2002black,rohwedder2011analysis}, which are prominent methods in self-consistent field theory~\cite{Ceniceros:2004uq,PhysRevLett.83.4317}. AA is also becoming popular in the numerical analysis community~\cite{walker2011anderson,toth2015convergence,evans2018proof,zhang2018globally,pollock2019anderson}. The literature on this subject is broad, so we only mention a few papers to show the variety of results obtained by AA.
In contrast to Picard iteration~\cite{picard1893application,butenko2014numerical}, which uses only one previous iterate, the method proceeds by linearly recombining a list of previous iterates in a way that approximately minimizes the linearized fixed-point residual. AA can be applied directly to accelerate fixed-point operators that arise naturally from solving partial differential equations (PDEs). The method was mainly used in optimization-free scenarios until the last few years. AA starts to show promising results in accelerating optimization algorithms~\cite{peng2018anderson,fu2019anderson,mai2019anderson,li2018fast} and machine learning~\cite{geist2018anderson}.

In this paper, we aim to combine the fast convergence of Newton-type methods with the low cost of only evaluating the gradient. We do this by applying an acceleration strategy introduced by D.G. Anderson~\cite{anderson1965iterative} to the steepest descent algorithm. We first reformulate the iterative formula as a fixed-point operator. In contrast to the classical gradient descent, AA produces a new iterate as a linear combination of several previous iterates. The linear coefficients are selected optimally to achieve the best reduction in the linearized fixed-point residual. As an acceleration strategy for the steepest descent method, AA can achieve competitive convergence speed with respect to methods like L-BFGS and nonlinear conjugate gradient descent (nCG) while reducing the computational cost of computing the exact Hessian or building an approximation of the (inverse) Hessian matrix. We illustrate the performance of these methods as the optimization algorithm for FWI and LSRTM in the numerical examples.

\section{Theory}
In this section, we first introduce the algorithmic details of AA for fixed-point problems and explain its similarities and differences with Picard iteration. Later, we review some essential background regarding FWI and LSRTM. Throughout the paper, we assume the forward model is an acoustic wave equation with a constant density.

\subsection{Anderson acceleration}
AA is an acceleration strategy introduced to improve the slow convergence of Picard iterations~\cite{anderson1965iterative}. We present the details of AA in Algorithm~\ref{alg:AA1}. The memory parameter $\cM$ determines the number of additional past iterates that need to be stored to compute the next iterate. For example, when $\cM=0$,  AA reduces to the Picard iteration as the $(k+1)$-th iterate $p_{k+1}$ only depends on $p_k$ and $p_{k+1} = G(p_k)$, where $G$ is the fixed-point operator. For a nonzero $\cM$, $p_{k+1}$ is a linear combination of the previous $\cM+1$ iterates, together with their evaluation by the fixed-point operator $G$; see~\Cref{eq:AAupdate0} for the detailed updating formula. If $\cM=+\infty$ and the damping parameter $\beta_k$ in~\Cref{eq:AAupdate0} is also chosen optimally at every iteration, AA is essentially equivalent to the Generalized Minimal Residual Method (GMRES) when $G$ is a linear fixed point operator and the fixed-point solution solves the square linear system $Ax=b$~\cite{toth2015convergence}. The damping parameter $\beta_k$, which could vary at different iteration $k$, controls the balance between the linear combination of the iterates $\{p_{k-\cM+i}\}_{i=0}^{\cM}$ and the linear combination of their evaluation by the operator $\{G(p_{k-\cM+i})\}_{i=0}^{\cM}$.

\begin{algorithm}[ht!]
  \caption{Anderson Acceleration \label{alg:AA1}}
  \begin{algorithmic}
       \STATE {\bf Input:} Given the initial guess $p_0$ and memory parameter $\cM \ge 1$. $G$ is the given fixed-point operator. Set $p_1 = G(p_0)$.
\FOR{k = 0, 1, 2,\ldots}
       \STATE {\bf Step 1:} Set $\cM_k = \min(\cM,k)$ and matrix $F_k = (f_{k-\cM_k},\ldots,f_k)$, where $f_i = G(p_i)-p_i$ is the fixed-point residual of the $i$-th iterate.
\STATE {\bf Step 2:} Find the optimal weights $\alpha^{(k)} = (\alpha^{(k)}_0,\ldots,\alpha^{(k)}_{\cM_k})^T$ by the optimization problem
\begin{equation}\label{eq:AAopt}
   \min\limits_{\sum_{i=0}^{\cM_k} \alpha_i^{(k)}=1} \| F_k  \alpha^{(k)} \|_*.
\end{equation}
\STATE {\bf Step 3:} Update the next iterate $p_{k+1}$
\bq \label{eq:AAupdate0}
p_{k+1} = (1-\beta_k) \sum_{i=0}^{\cM_k} \alpha_i^{(k)} p_{k-\cM_k+i}
 + \beta_k \sum_{i=0}^{\cM_k} \alpha_i^{(k)} G(p_{k-\cM_k+i}).
\eq
\ENDFOR
    \end{algorithmic}
\end{algorithm}

At iteration $k$, the coefficient vector $\alpha^{(k)}=(\alpha^{(k)}_0,\ldots,\alpha^{(k)}_{\cM_k})^T$ is determined by minimizing the sum of the weighted fixed-point residuals. The sum of all the coefficients must total one so that the fixed-point solution $p^*$ is preserved under the updating formula~\eqref{eq:AAupdate0} of AA. The Picard iteration fits into the updating formula~\eqref{eq:AAupdate0} with the weighting vector $(0,0,\dots,0,1)^T$ for every $k$. Since the weighting vector for AA is obtained from the optimization problem~\eqref{eq:AAopt}, AA is always at least as good as Picard iteration as
\begin{equation}
    \bigg\| \sum_{i=0}^{\cM_k} \alpha_i^{(k)} f_{k-\cM_k+i}\bigg\|_* \leq \| f_k\|_*,
\end{equation}
where $f_i = G(p_i) - p_i$ is the fixed-point residual of the $i$-th iteration.

Through a change of variable, one can remove the constraints in~\eqref{eq:AAopt} to simplify the optimization step. Consider a new vector $\gamma^{(k)}  = (\gamma_0^{(k)} ,\dots,\gamma_{\cM_k-1}^{(k)} )^T$ defined by the optimal parameter  $\alpha^{(k)}$ where
\begin{equation}
\gamma_i^{(k)} = \alpha_0^{(k)} + \dots + \alpha_i^{(k)},\ 0\leq i \leq \cM_k-1. 
\end{equation}
Consider the matrix $A_k$ given by
\begin{equation}~\label{eq:Ak}
A_k = ( f_{k-\cM_k+1} -  f_{k-\cM_k},\dots,  f_{k} - f_{k-1}),
\end{equation}
whose column vectors are the differences in the fixed-point residual between two consecutive iterations. The optimization step~\eqref{eq:AAopt} is equivalent to the following unconstrained optimization problem
\begin{equation}\label{eq:AAopt2}
\gamma^{(k)} = \argmin\limits_\gamma \| A_k \gamma -f_k\|_*.
\end{equation}

There are several variants regarding the choice of the norm $\|\cdot \|_*$ in the optimization step~\eqref{eq:AAopt}. For example, one can use the $\ell^1$, $\ell^2$, or the $\ell^\infty$ norm as the objective function. Alternatively, a weighted $\ell^2$ norm may improve the conditioning of the fixed-point operator or enforce the spectral bias towards certain modes of the solution~\cite{yang2020anderson}. The optimal weights may not be the same among different choices of the objective function, and the cost of solving the corresponding optimization problem can also be radically different. For example, linear programming is required to solve the optimization under the $\ell^1$ and the $\ell^\infty$ norms. If we stick with the $\ell^2$ norm, $ \gamma^{(k)}$ is then the least-squares solution to the following linear system:
\begin{equation}\label{eq:solveGamma}
A_k  \gamma^{(k)}  =  f_k,
\end{equation} 
where $ f_k = G(p_k) - p_k$ is the fixed-point residual at iteration $k$. If $\beta_k = 1$ for any $k$, then we can rewrite the updating formula in terms of $\gamma^{(k)}$ as follows:
\begin{equation} \label{eq:AA2_update}
p_{k+1}  = G(p_k) -  \sum_{i=0}^{\cM_k-1}  \gamma_i^{(k)} \left[G(p_{k-\cM_k+i+1}) - G(p_{k-\cM_k+i}) \right]. 
\end{equation}
It is computationally efficient to implement AA based on~\eqref{eq:AA2_update}. The size of $A_k$ is $n$ by $\cM_k$, where $n$ is the dimension of the parameter and $\cM_k = \min(\cM,k) \leq \cM$. The memory parameter $\cM$ of AA is often chosen to be small. A rank-updated QR factorization can further reduce the cost of solving~\eqref{eq:solveGamma}~\cite[ Section 12.5.1]{golub2012matrix}.

\subsection{FWI and LSRTM}
Seismic inversion aims to obtain an estimate of the distribution of material properties in the underground. They are large-scale inverse problems that we treat as constrained optimization problems based on the deterministic approach of solving inverse problems.

FWI is a nonlinear inverse technique that
utilizes the entire wavefield information to estimate
the earth's properties. Without loss of generality, the PDE constraint of FWI is the following acoustic wave equation with zero initial condition and non-reflecting boundary conditions.
\begin{equation}\label{eq:FWD}
     \left\{
     \begin{array}{rl}
     & m(\mathbf{x})\dfrac{\partial^2 u(\mathbf{x},t)}{\partial t^2}- \Laplace u(\mathbf{x},t) = s(\mathbf{x},t),\\
    & u(\mathbf{x}, 0 ) = 0,                \\
    & \dfrac{\partial u}{\partial t}(\mathbf{x}, 0 ) = 0 . 
     \end{array} \right.
\end{equation}
We set the model parameter $m(\mathbf{x}) = \dfrac{1}{c(\mathbf{x})^2}$, where $c(\mathbf{x})$ is the wave velocity, $u(\mathbf{x},t)$ is the forward wavefield, $s(\mathbf{x},t)$ is the wave source. The velocity parameter $m$ is often the target of reconstruction. Equation~\eqref{eq:FWD} is a linear PDE but defines a nonlinear operator $\mathcal F$ that maps $m(\mathbf{x})$ to $u(\mathbf{x},t)$. In FWI, we translate the inverse problem of finding the model parameter based on the observable seismic data to a constrained optimization problem:
\begin{equation}\label{eq:FWI}
    m^{\ast} = \argmin_{m} J(m),\quad J(m) = \dfrac{1}{2}||f(m)-g||_2^2.
\end{equation}
The least-squares norm is commonly used as the objective function $J$ to calculate the misfit between the synthetic data $f(m) = R\mathcal F(m)$ and the observed data $g$. Here, $R$ is the projection operator that extracts the wavefield $u$ at the receiver locations. There are other choices of objective functions to mitigate the cycle-skipping issues of FWI~\cite{yang2018application}. 

LSRTM is a new migration method designed to improve the image quality generated by reverse-time migration (RTM). It is formulated as a linear inverse problem based on the Born approximation, a first-order linearization of the map $\mathcal F$~\cite{hudson1981use}. From now on, we denote the forward operator of LSRTM, i.e., the Born modeling, as $\mathcal L = \dfrac{\delta \mathcal F}{\delta m}$, the functional derivative of $\mathcal F$ with respect to $m$. The linear operator $\mathcal L$ maps a small perturbation in the velocity $m_r$ to the scattering wavefield $u_r$:
\begin{equation}\label{eq:Born}
     \left\{
     \begin{array}{rl}
     & m_0\dfrac{\partial^2 u_r(\mathbf{x},t)}{\partial t^2} - \Laplace u_r(\mathbf{x},t) = -m_r \dfrac{\partial^2 u_0(\mathbf{x},t)}{\partial t^2}, \\
    & u_r(\mathbf{x}, 0 ) = 0,                \\
    & \dfrac{\partial u_r}{\partial t}(\mathbf{x}, 0 ) = 0. 
     \end{array} \right.
\end{equation}
Here, $m_0$ is the given background velocity and the background wavefield $u_0 = \mathcal F(m_0)$. We seek the reflectivity model by minimizing the least-squares error between the observed data $d_r$ and the predicted scattering wavefield $L m_r = R\cL m_r$,
\begin{equation}~\label{eq:LSRTM}
m_r^\ast = \argmin_{m_r} J(m_r),\quad J(m_r)= \dfrac{1}{2} ||L m_r - d_r||_2^2.
\end{equation}

To solve for $m^\ast$ in~\Cref{eq:FWI} and $m_r^\ast$ in~\Cref{eq:LSRTM}, optimization algorithms heavily rely on the gradient and the Hessian information of the objective function $J$. In seismic inversions, one can obtain the gradient of a parameter by solving the forward equation and the adjoint equation once, based on the adjoint-state method~\cite{Plessix}. The adjoint equation for both FWI and LSRTM is the following
\begin{equation} \label{eq:FWI_adj}
     \left\{
     \begin{array}{rl}
     & m\dfrac{\partial^2 v(\mathbf{x},t)}{\partial t^2}- \Laplace v(\mathbf{x},t)  = -R^*\dfrac{\partial J}{\partial f},\\
    & v(\mathbf{x}, T ) = 0,                \\
    & \dfrac{\partial v}{\partial t}(\mathbf{x}, T ) = 0. 
     \end{array} \right.
\end{equation}
For LSRTM, the $m$ in~\Cref{eq:FWI_adj} is the background velocity $m_0$. The term $\dfrac{\partial J}{\partial f}$ is the Fr\'{e}chet derivative of the objective function with respect to the synthetic data $f$, also known as the adjoint source. If $J$ is the least-squares norm, $\dfrac{\partial J}{\partial f} = f-g$, simply the data residual. The functional derivative of the objective function $J$ with respect to the model parameter $m$ is
\begin{equation}~\label{eq:adj_grad1}
\frac{\partial J}{\partial m}  =- \int_0^T \frac{\partial^2 u(\mathbf{x},t)}{\partial t^2} v(\mathbf{x},t)dt,
\end{equation}  
where $u$ and $v$ are the forward and adjoint wavefields, respectively. For FWI, $u$ is the solution to the acoustic wave equation~\eqref{eq:FWD}, while for LSRTM, $u$ is the solution to linearized wave equation~\eqref{eq:Born}. An outstanding advantage of the adjoint-state method is that the number of wave simulations to compute the gradient is independent of the size of $m$. The model parameter can then be updated by a gradient-based optimization algorithm iteratively till meeting the stopping criteria. The Hessian matrix can also be computed based on the adjoint-state method if it is needed for optimization, uncertainty quantification or resolution analysis.


\section{Method}
Within the framework of iterative methods, we treat the gradient descent algorithm as a fixed-point iteration, where the fixed-point solution is the optimal model parameter. Consider an objective function $J(p)$ for the unknown parameter $p$. If we choose to minimize $J(p)$ by the steepest descent algorithm, then the $(k+1)$-th iterate $p_{k+1}$ is obtained by the $k$-th iterate $p_k$ and the gradient vector for $p_k$. The step size is properly chosen to guarantee a sufficient decrease in the objective function. Without loss of generality, we fix the step size as a small positive constant $\eta$, and obtain the updating formula by steepest descent:
\begin{equation} \label{eq:GD}
p_{k+1} = p_k - \eta \frac{\partial J}{\partial p}\bigg|_{p=p_k} = G(p_k).
\end{equation}
Since the right-hand side of~\eqref{eq:GD} only depends on $p_k$, one can regard the updating formula as a fixed-point operator $G$ applied to $p_k$. Equation~\eqref{eq:GD} can be considered as the Picard iteration for $G$. The fixed-point solution $p^*$ should satisfy
\begin{equation}
    p^* = G(p^*)\quad \Longleftrightarrow \quad \frac{\partial J}{\partial p}\bigg|_{p=p^*} =0.
\end{equation} 
Thus, $p$ is the fixed-point solution of $G$ if and only if the gradient of the objective function $J$ is zero at $p$.

Seismic inverse problems are often ill-posed and suffer from local minima trapping. Typically, the zero-gradient condition is far from enough to guarantee the optimality, especially for FWI. There has been extensive literature on tackling the nonconvexity~\cite{engquist2020optimal,symes2020full}. In this paper, we focus on accelerating the convergence and not addressing the cycle-skipping issues, which is another important research topic by itself. Thus, we assume the initial guess $p_0$ in this paper is good enough that the optimization problem does not suffer from local minima trapping. 

Given the fixed-point operator $G$ defined by the steepest descent algorithm~\eqref{eq:GD}, we aim to accelerate the convergence by applying the Anderson acceleration. First, we rewrite~\Cref{eq:Ak} as follows:
\begin{equation}
    f_k = G(p_k) - p_k = -\eta \mathcal{G}_k,
\end{equation}
\begin{equation}
A_k= -\eta (\mathcal{G}_{k-\cM_k+1}-\mathcal{G}_{k-\cM_k},\dots,\mathcal{G}_k-\mathcal{G}_{k-1}),
\end{equation}
where $\mathcal{G}_k = \frac{\partial J}{\partial p}\big|_{p_k}$ is the gradient vector for the iterate $p_k$. We recall that the core of AA is to solve a linear system~\eqref{eq:solveGamma} for $\gamma^{(k)}$ which gives the optimal coefficients $\alpha^{(k)}$ by a change of variable. Applying AA to accelerate gradient descent, we remark that the main components of the linear system~\eqref{eq:solveGamma} are constructed only by the gradients $\{\mathcal G_i\}_{i=k-\cM_k}^{k}$ of the optimization. Thus, the memory requirements of AA is the same as the L-BFGS algorithm.

For typical fixed-point problems, the fixed-point residual is the indicator of convergence. That is, we judge the convergence by comparing the norm of $G(p)-p$. As a unique feature of our application,  the fixed-point operator comes from an optimization problem. Thus, the objective function can also be utilized in AA to improve the convergence further. The traditional AA described in~\Cref{alg:AA1} does not have such a step related to the objective function. Therefore, 
by combining both the fixed-point residual and the objective function, we describe a new workflow of AA for seismic inversion in~\Cref{alg:AA2}. The final $p_{k+1}$ is a linear combination of the output by the gradient descent, $\bar{p}_{k+1}$, and the optimized new iterate by AA, $\widetilde{p}_{k+1}$. A backtracking line search following the Wolfe condition is applied to determine the weighting between $\bar{p}_{k+1}$ and $\widetilde{p}_{k+1}$ such that we can achieve the best decrease in the objective function $J$; see~\Cref{alg:AA2} for more details.

\begin{algorithm}[ht!]
  \caption{$\ell^2$-based AA for gradient descent\label{alg:AA2}}
  \begin{algorithmic}
  \STATE {\bf Input:} Given the initial model $p_0$, the memory parameter $\cM \ge 1$, and the fixed-point operator $G$ based on the gradient descent update~\eqref{eq:GD}. Set $p_1 = G(p_0)$.
\FOR{$k = 1, 2, \ldots$ until convergence or maximum iteration}
       \STATE {\bf Step 1: }Set $\bar{p}_{k+1} = G(p_k) = p_k - \eta \mathcal G_k$. Update $A_k$ and $f_k$ using the gradient vectors $\{\mathcal G_i\}_{i=k-\cM_k}^{k}$ where $\cM_k = \min(\cM,k)$.
\STATE {\bf Step 2:} Find the least-squares solution $\gamma^{(k)}$ to the linear system~\eqref{eq:solveGamma} by using the low-rank QR update.
\STATE {\bf Step 3: }Compute the new iterate  $\widetilde{p}_{k+1}$ following~\Cref{eq:AA2_update}.
\STATE {\bf Step 4: }Apply the backtracking line search for $\lambda$ such that $J(\lambda \widetilde{p}_{k+1} +  ( 1-\lambda) \bar{p}_{k+1})$ has a sufficient decrease compared to $J(p_k)$.
\STATE  {\bf Step 5: }Set the new iterate as
\begin{equation}
    p_{k+1} =  \lambda  \widetilde{p}_{k+1} +  ( 1- \lambda ) \bar{p}_{k+1}.
\end{equation}
\ENDFOR
    \end{algorithmic}
\end{algorithm}


\section{Numerical Example}
In this section, we present several inversion tests for both FWI and LSRTM and compare the performance of AA with other optimization methods such as L-BFGS and steepest descent.

\subsection{Full waveform inversion}
\begin{figure}
\begin{center}
\subfloat[True velocity]{\includegraphics[width=0.5\textwidth]{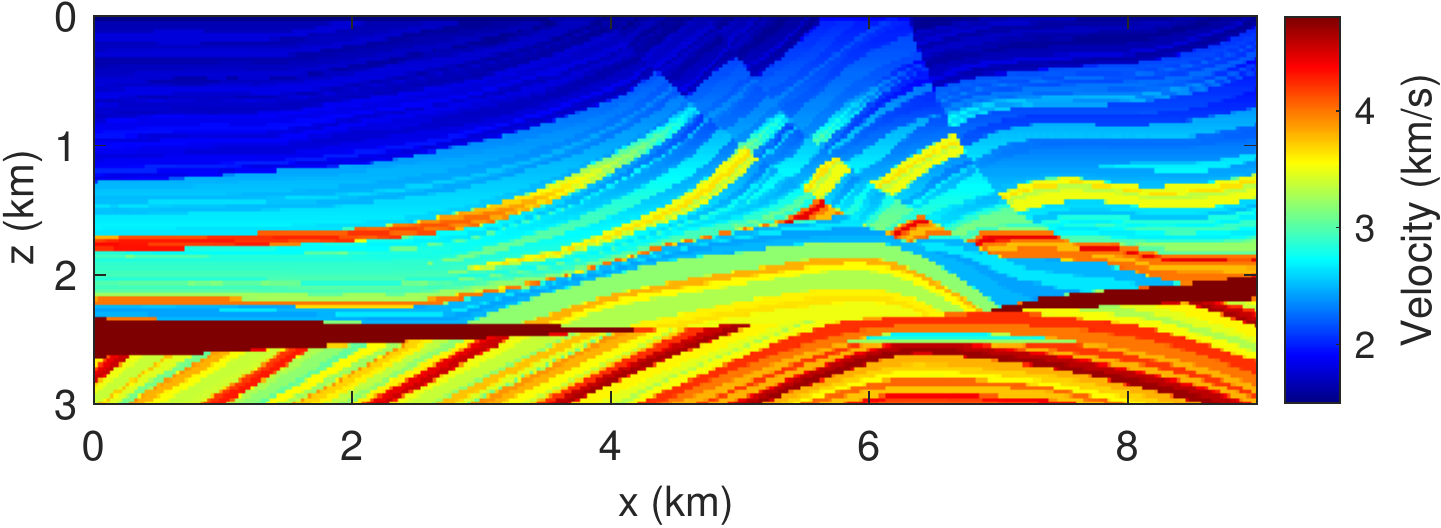}\label{fig:Marm-true}}
\subfloat[Initial velocity]{\includegraphics[width=0.5\textwidth]{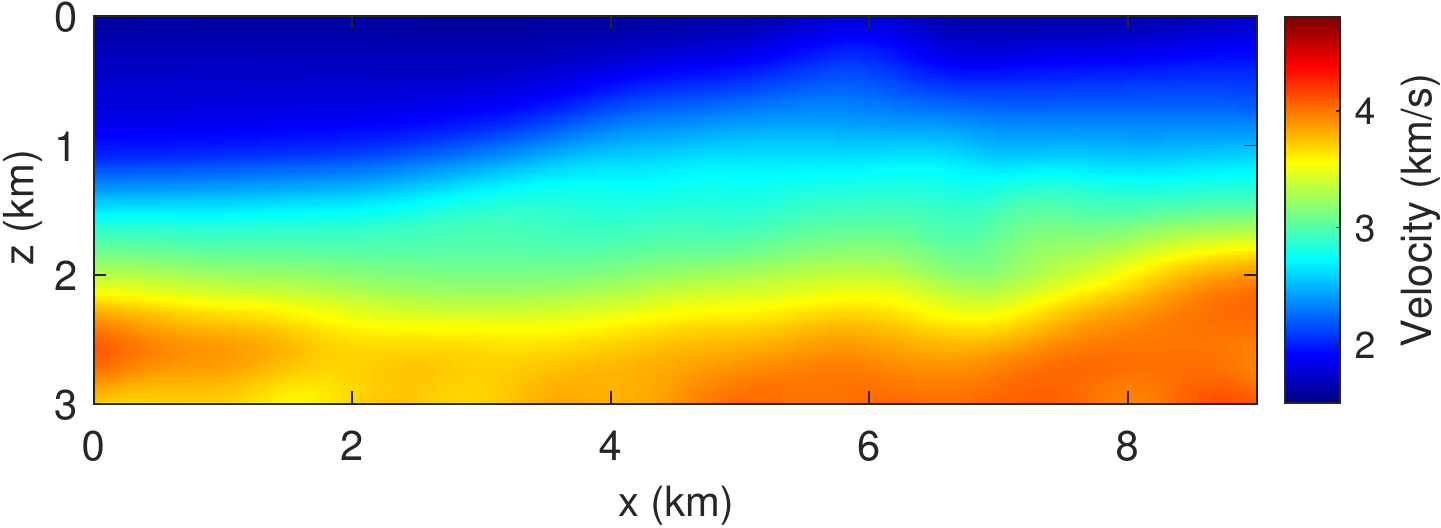}\label{fig:Marm-v0}}
\caption{Marmousi true and initial velocity models for the FWI tests.}~\label{fig:FWI-true}
\end{center}
\end{figure}

We aim to reconstruct the  Marmousi velocity model that is 3 km in depth and 9 km in width (\Cref{fig:Marm-true}) from a smoothed initial guess as shown in~\Cref{fig:Marm-v0}. There are $11$ equally spaced sources at $150$ m below the air-water interface. The source is a Ricker wavelet centered at 15 Hz, and 4 seconds is the total recording time. There is no cycle skipping with the chosen initial model. The most time-consuming component of FWI is seismic modeling, which is essential for gradient calculation. Thus, instead of counting the number of iterations, we use the number of gradient evaluations as a measure of performance.

\begin{figure}
\begin{center}
\subfloat[AA ($\cM=20$)]{ \includegraphics[width=0.5\textwidth]{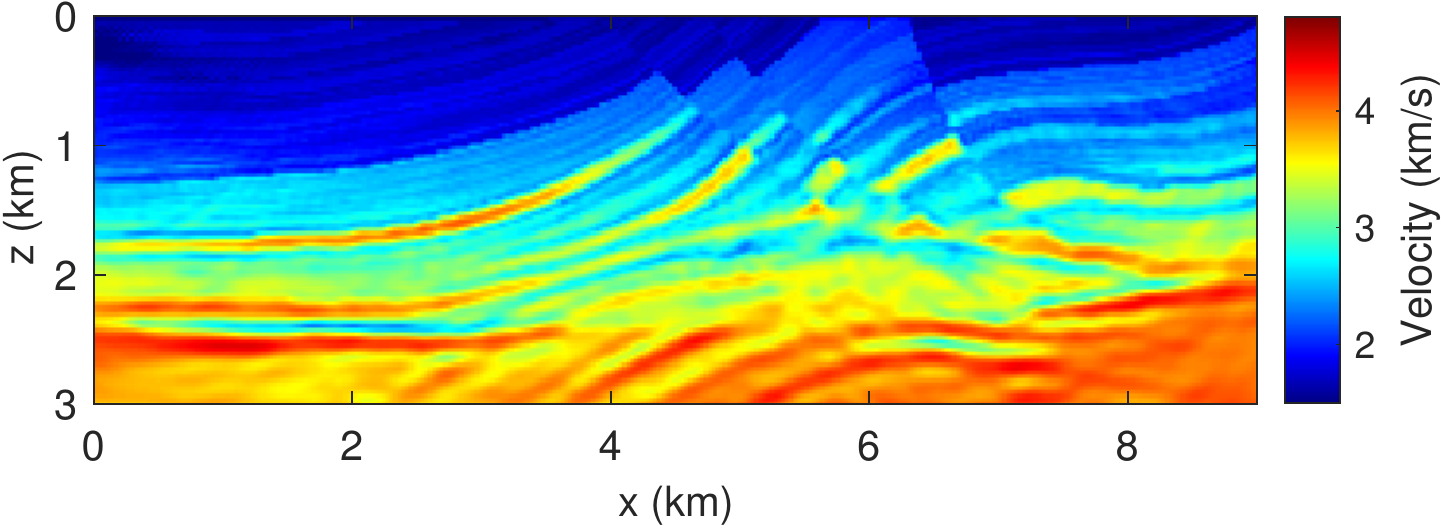}\label{fig:inv-aa}}
\subfloat[L-BFGS ($\cM=20$)]{
\includegraphics[width=0.5\textwidth]{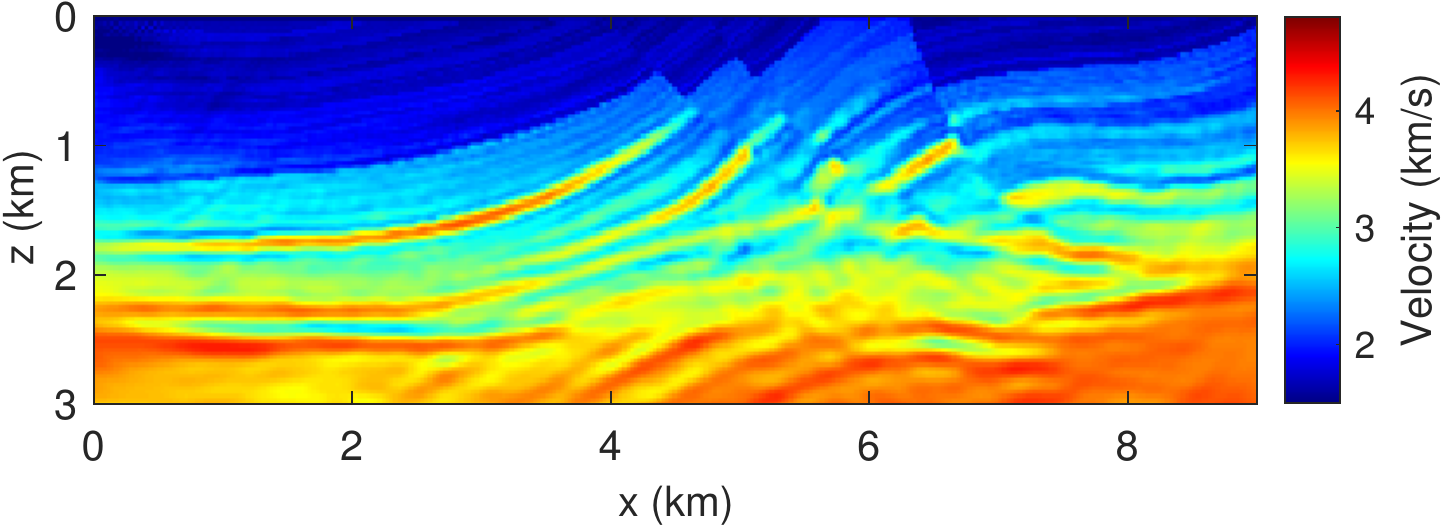}\label{fig:inv-bfgs}}\\
\subfloat[Nonlinear CG]{\includegraphics[width=0.5\textwidth]{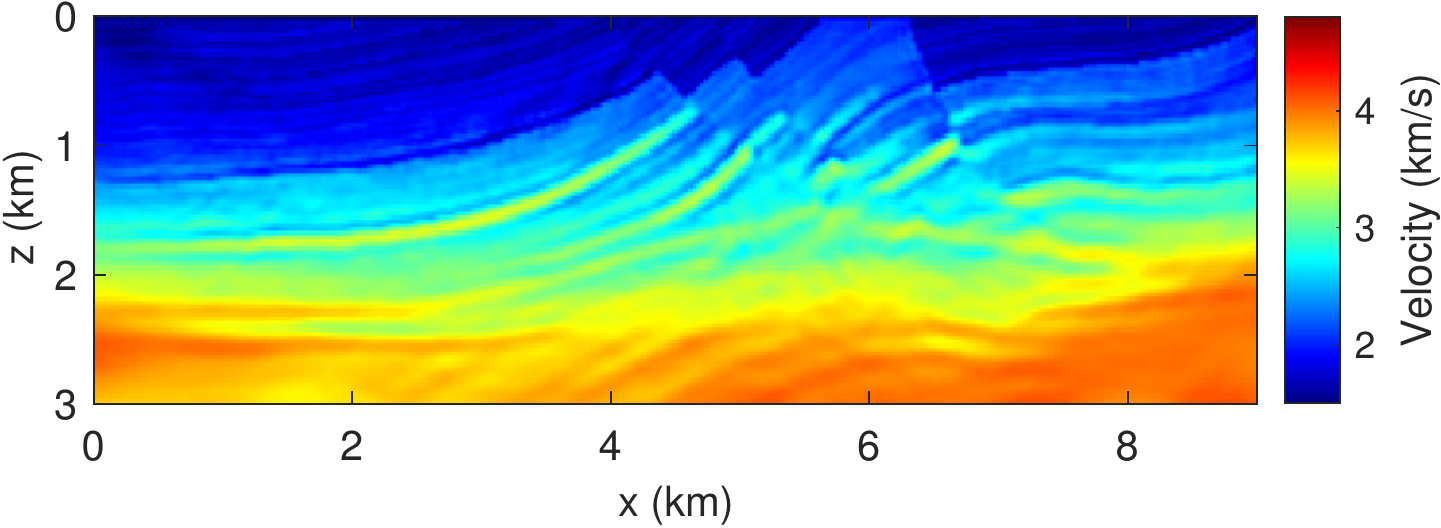}\label{fig:inv-ncg}}
\subfloat[Steepest descent]{\includegraphics[width=0.5\textwidth]{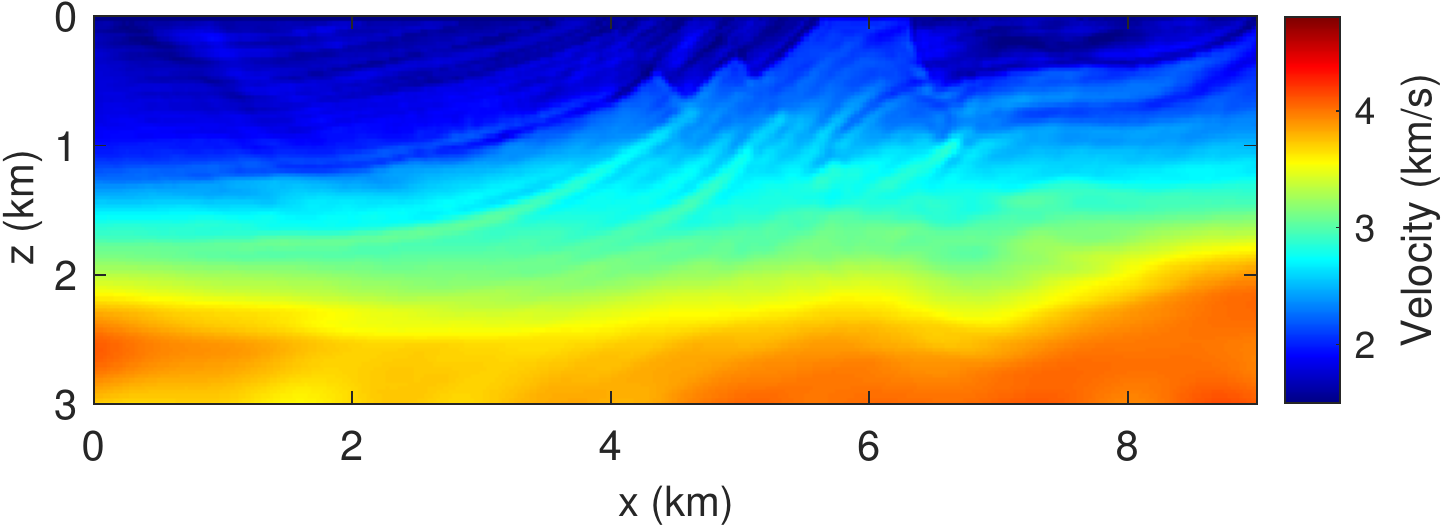}\label{fig:inv-sd}}
\caption{FWI results using AA (top left), L-BFGS (top right), nonlinear CG (bottom left) and steepest descent (bottom right) after $1000$ gradient evaluations.}~\label{fig:FWI-inv}
\end{center}
\end{figure}

\begin{figure} 
\begin{center}
\includegraphics[width=0.48\textwidth]{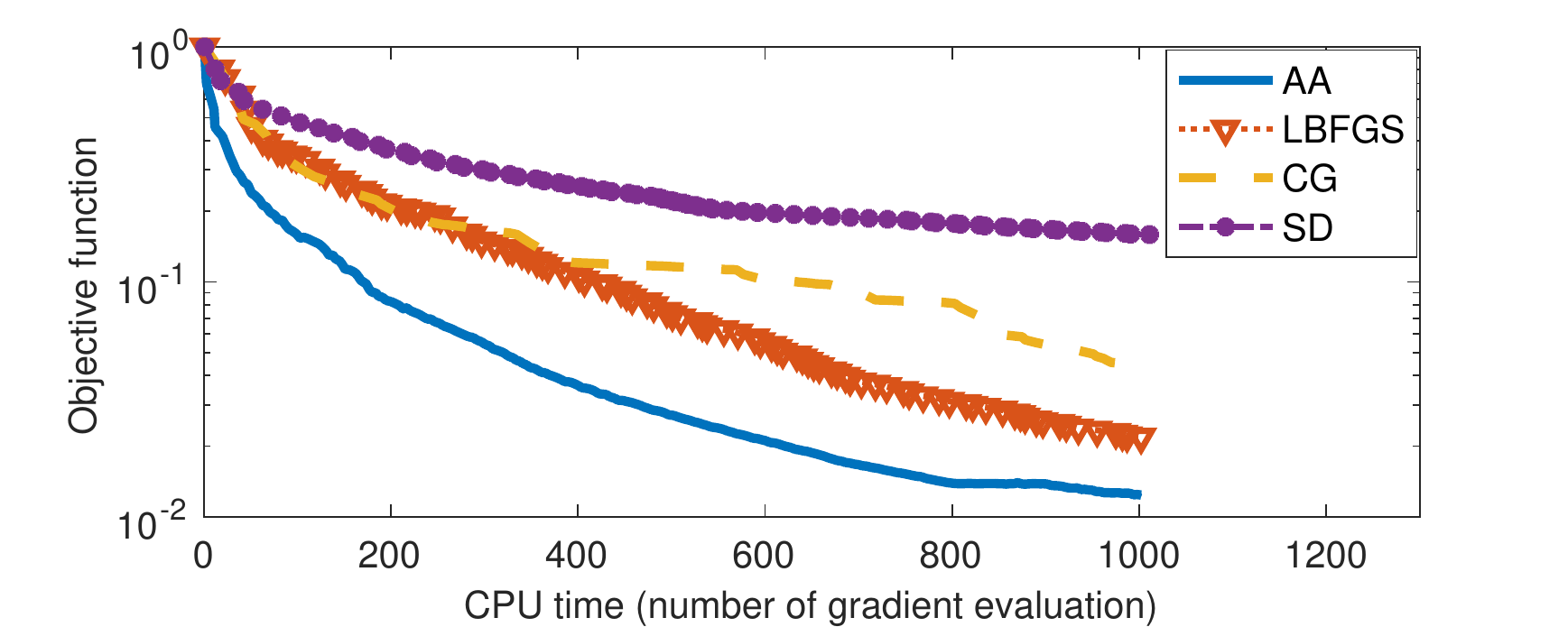}
\includegraphics[width=0.48\textwidth]{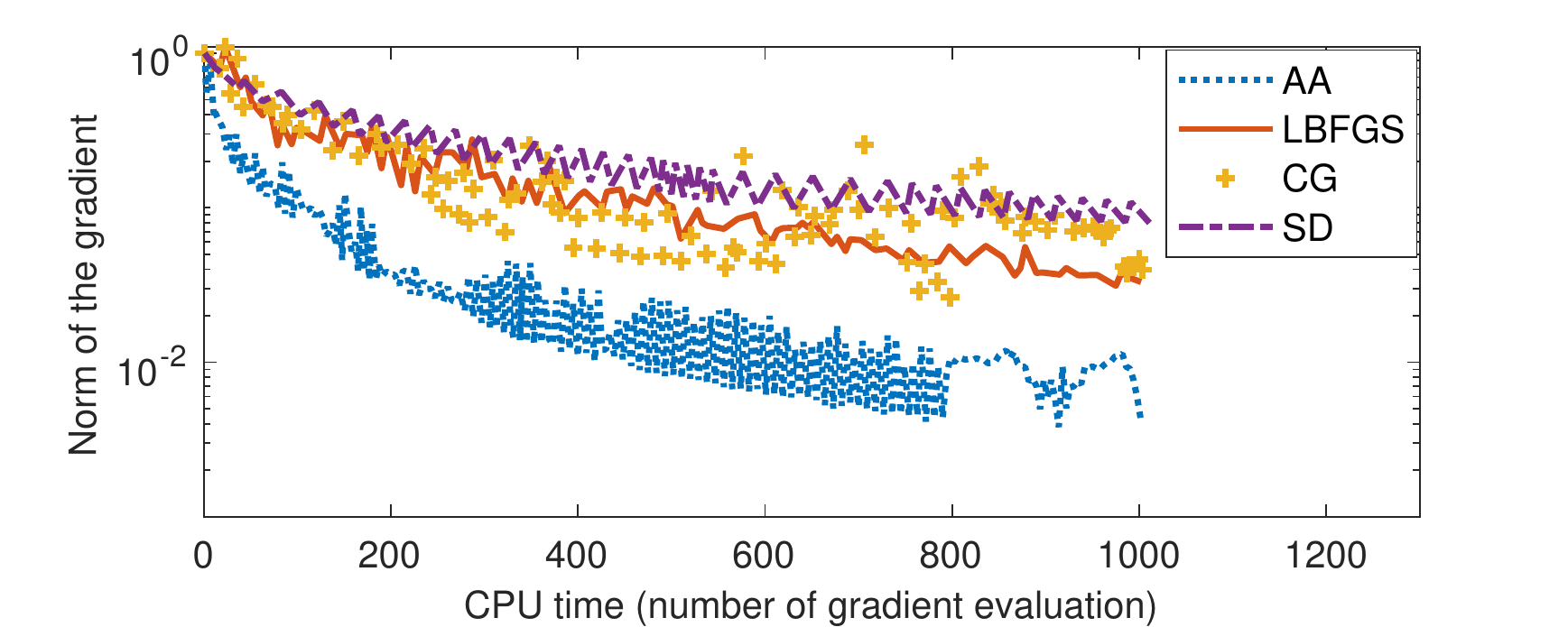}
\end{center}
\caption{FWI convergence history in terms of computational time (measured by the number of FWI gradient evaluations).}~\label{fig:FWI-curve}
\end{figure}

After 1000 gradient evaluations, FWI using AA ($\cM=20$), L-BFGS ($\cM=20$), nCG, and steepest descent are shown in~\Cref{fig:FWI-inv}. The same backtracking line search following the Armijo rule, and the curvature condition is applied to all methods. 
Results by AA and L-BFGS illustrate better resolution than the one by nCG. The steepest descent method converges slowly. The convergence history for both the $\ell^2$ objective function and the norm of the gradient is shown in~\Cref{fig:FWI-curve}. In both plots, AA demonstrates a faster convergence rate than both L-BFGS and nCG. Known as quasi-Newton methods, L-BFGS and CG converge in fewer iterations than AA. However, more gradient evaluations are spent on the backtracking line search, which increases the overall CPU time. The drastic improvement in the convergence rate between AA and the steepest descent method shows the benefits of this simple strategy by linearly recombining previous iterates. Considering the low cost of implementation, AA can be an attractive optimization technique for FWI. More analysis between AA and L-BFGS is presented in the next section.

\begin{figure}
\begin{center}
\subfloat[Data comparison]{ \includegraphics[width=0.48\textwidth]{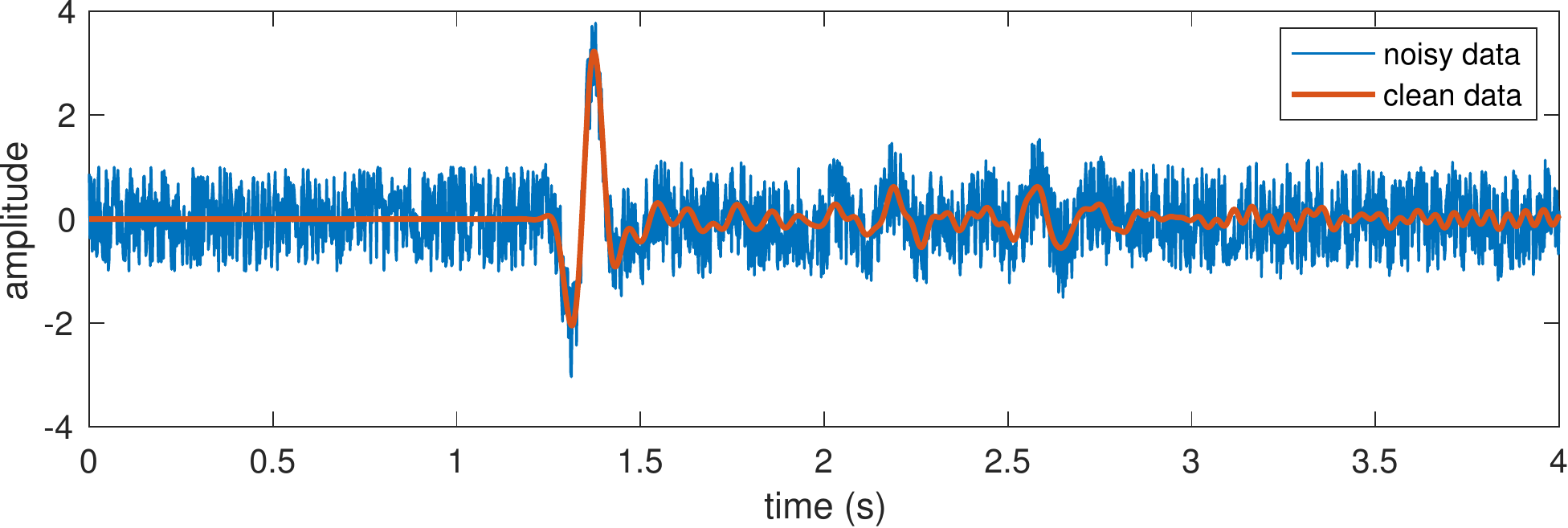}\label{fig:noise-data}}
\subfloat[AA ($\cM=20$)]{
\includegraphics[width=0.5\textwidth]{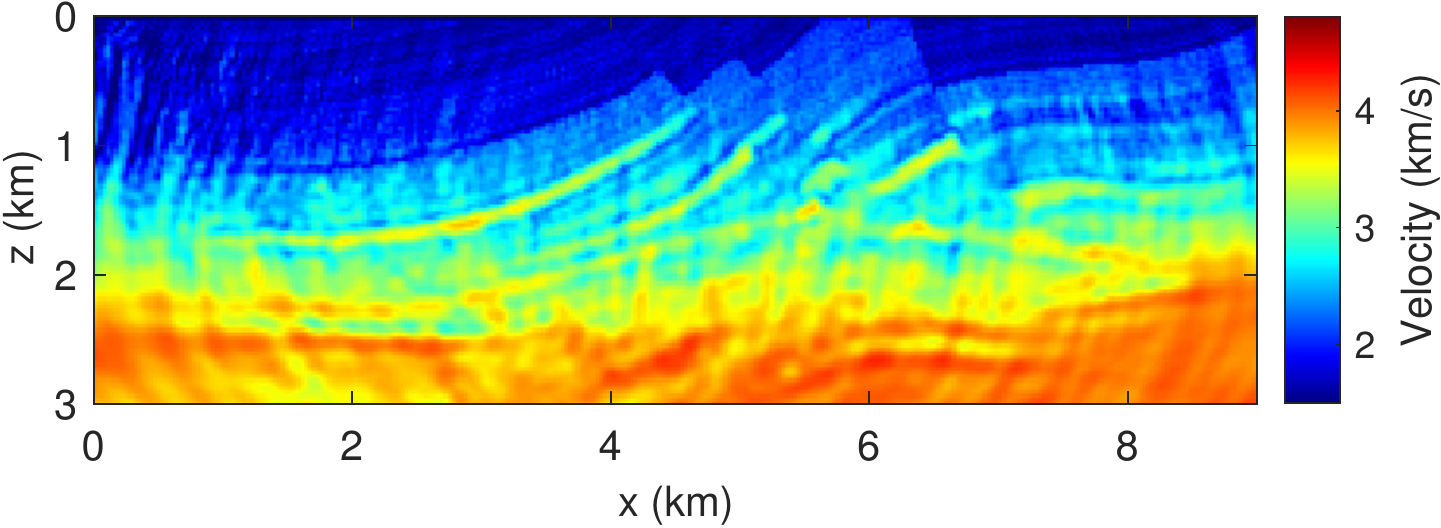}\label{fig:noise-aa}}\\
\subfloat[L-BFGS ($\cM=20$)]{\includegraphics[width=0.5\textwidth]{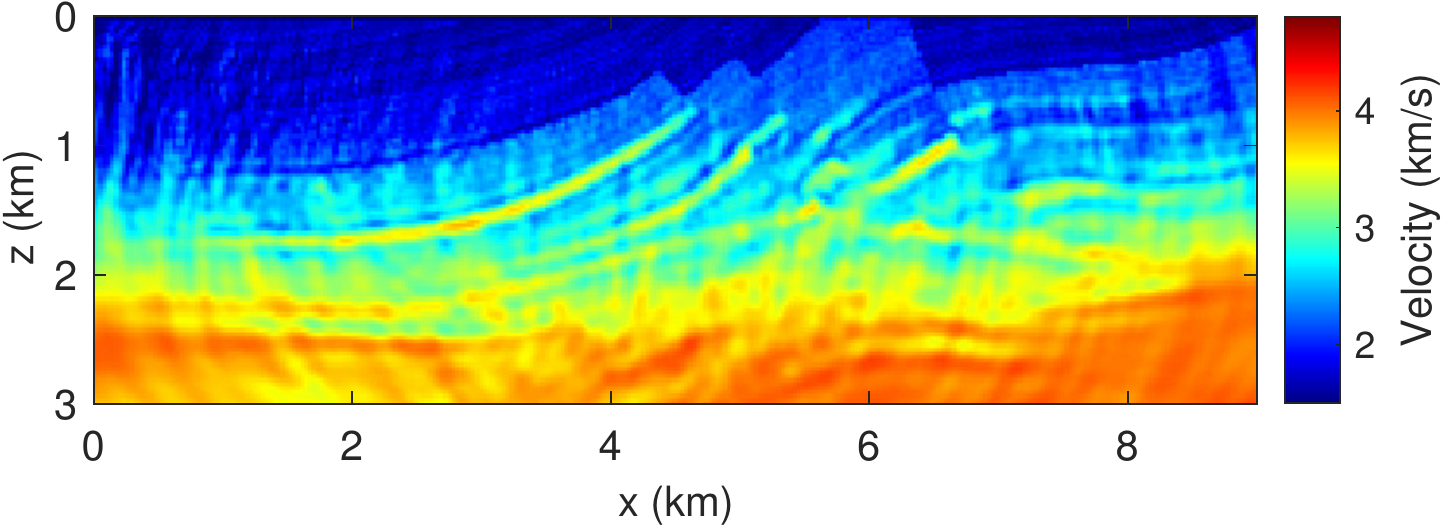}\label{fig:noise-bfgs}}
\subfloat[Steepest descent]{\includegraphics[width=0.5\textwidth]{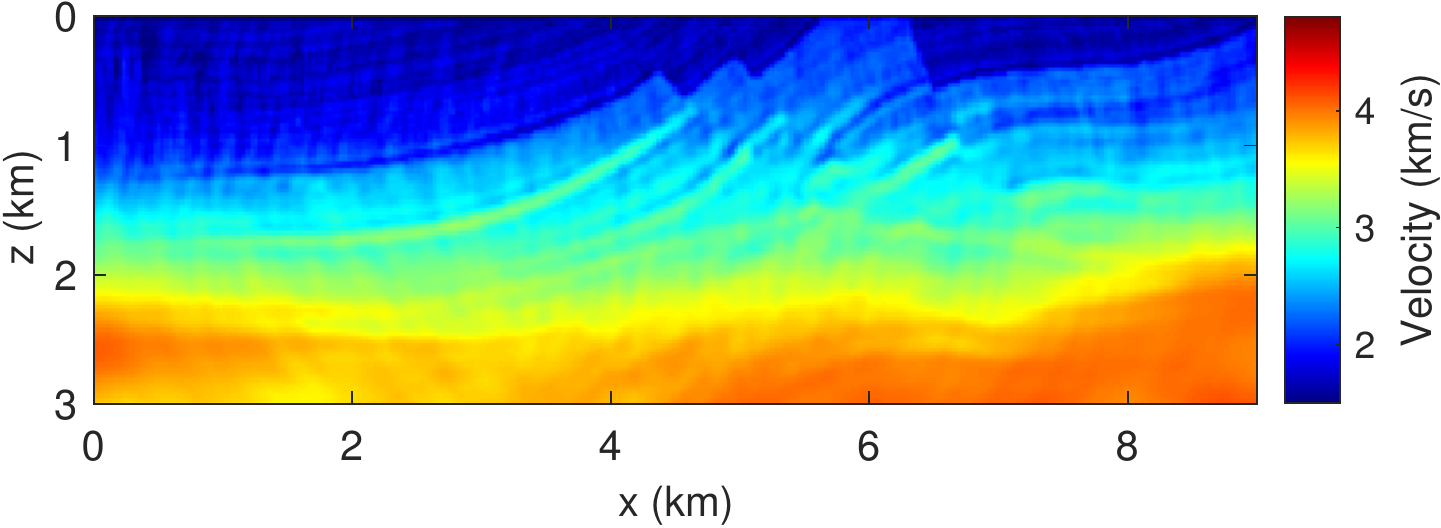}\label{fig:noise-sd}}
\caption{FWI with noise: (a) the comparison between one trace of noisy and clean signals from the observed data; (b) the inversion using AA, (c) the L-BFGS algorithm and (d) the steepest descent method after $1000$ gradient evaluations.}~\label{fig:FWI-noise-inv}
\end{center}
\end{figure}

\subsection{Inversion with noise}
We present another FWI example. Unlike the previous set of tests, we add mean-zero noise to the observed data to make the inversion test more representative of results expected on real data. The noise follows a uniform distribution, and the signal-to-noise ratio (SNR) is $0.55$ dB. We plot one trace of the clean data and the noisy data in~\Cref{fig:noise-data} for illustration. All the other settings remain the same as the noise-free example.~\Cref{fig:noise-aa} is the inversion result using AA to accelerate the steepest descent algorithm while \Cref{fig:noise-sd} is the inversion result without acceleration. We also perform a test for the L-BFGS algorithm; see~\Cref{fig:noise-bfgs}. All experiments are stopped after $1000$ gradient calculations as before. Compared with the noise-free results in the previous section, one can observe artificial oscillatory features in all the images resulting from noise overfitting. However, the noise footprints are equally strong for the inversion using L-BFGS and the one using AA. Although the reconstruction by the steepest descent method seems to be less noisy, it also recovers fewer features of the true Marmousi model. It is due to its slower convergence compared with AA and L-BFGS. Typically, one can expect that the artifacts in the reconstructed model to be proportional to the noise SNR. To mitigate the noise effects, one can change the objective function from the $\ell^2$ norm to the $W_2$ metric~\cite{yang2018application}, which is proved to be more robust with respect to noise. One can also add regularization terms to the objective function, which is a common strategy to improve the stability of the inverse problem.


\subsection{Least-squares reverse-time migration}
\begin{figure} 
\begin{center}
\subfloat[Background velocity for RTM/LSRTM]{ \includegraphics[width=0.7\textwidth]{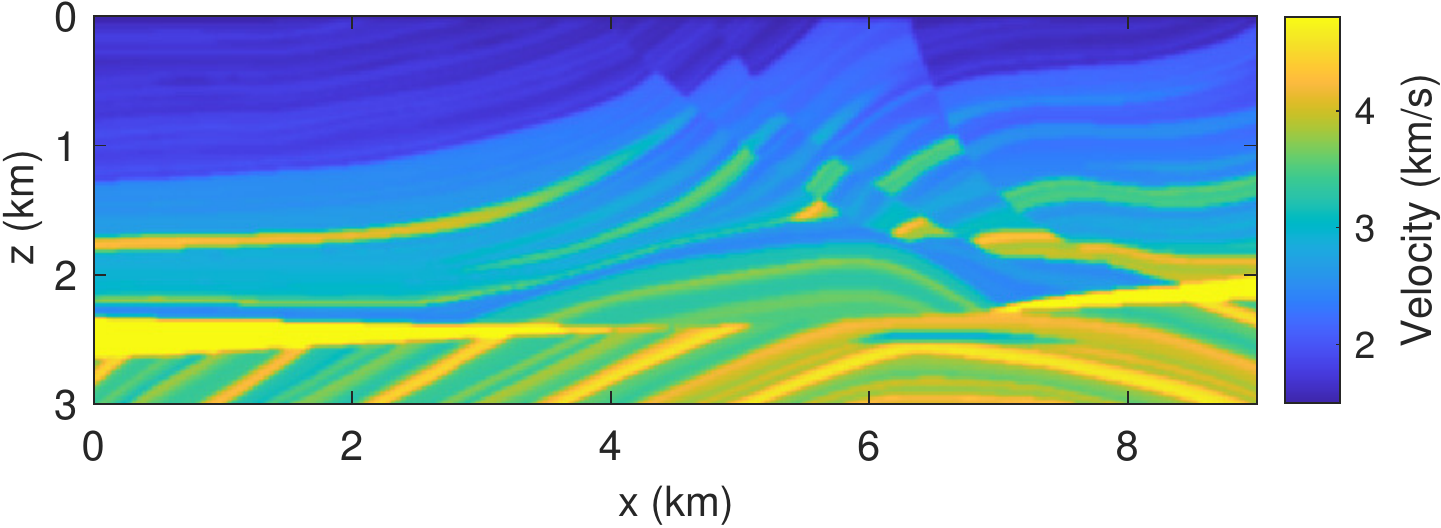}~\label{fig:Mig-smooth}}\\
\subfloat[True reflectivity]{ \includegraphics[width=0.7\textwidth]{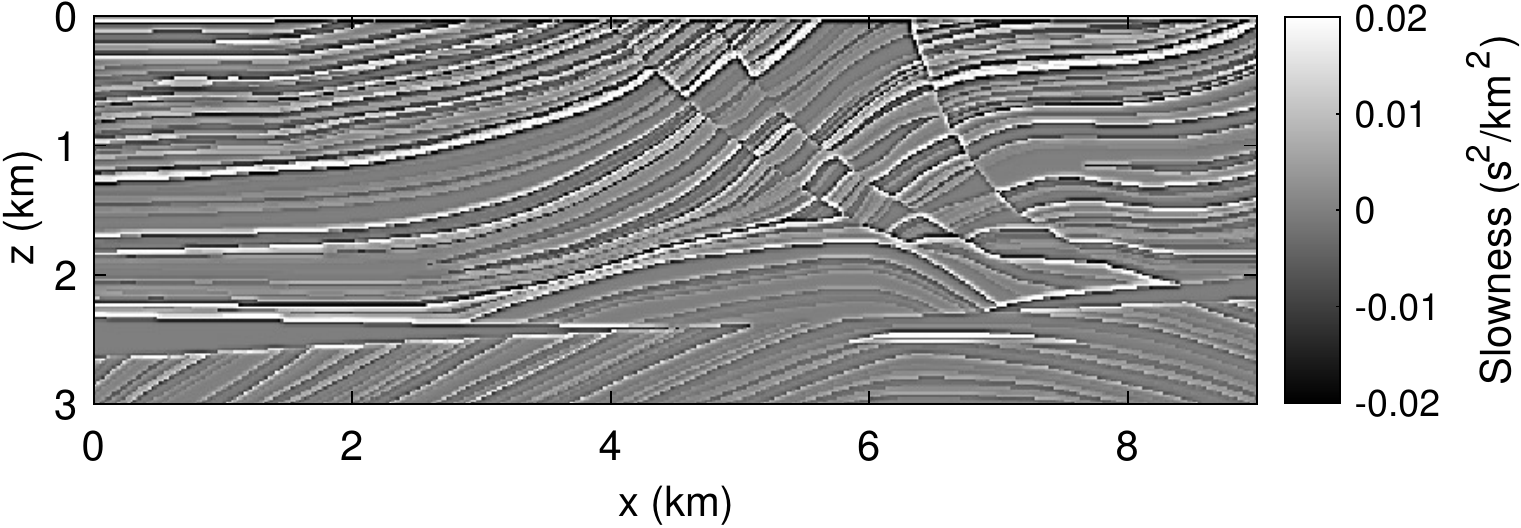}~\label{fig:Mig-true}}\\
\subfloat[RTM result]{ \includegraphics[width=0.7\textwidth]{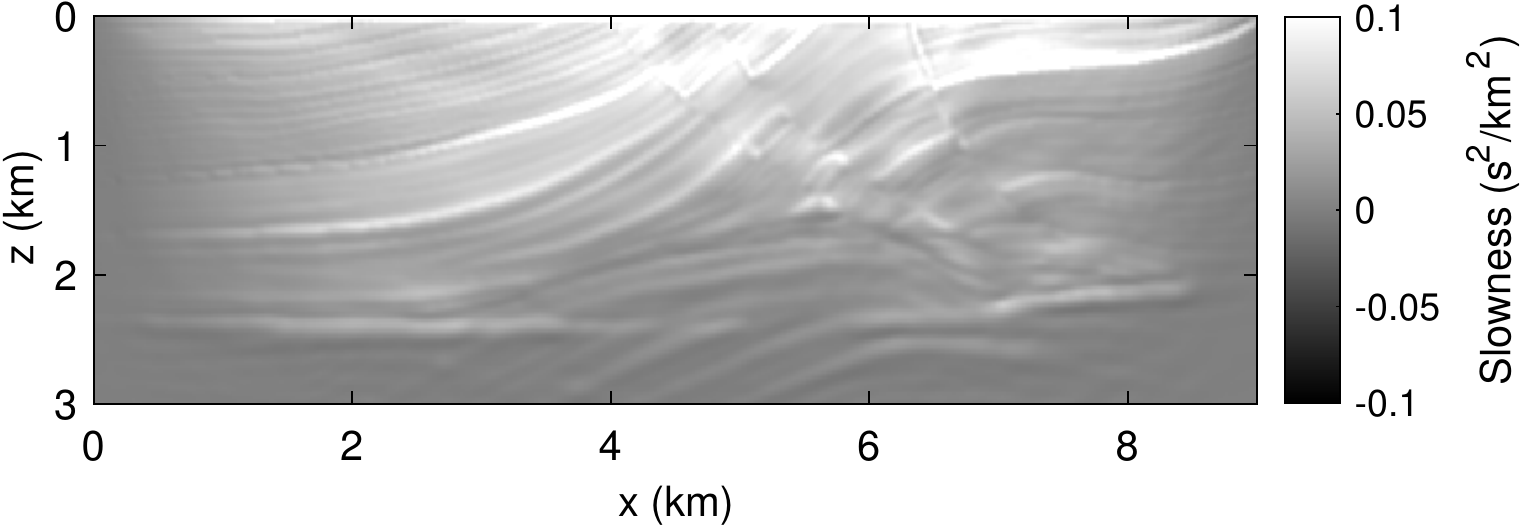}\label{fig:Mig-RTM}}\\
\subfloat[RTM result after Laplacian filtering]{ \includegraphics[width=0.7\textwidth]{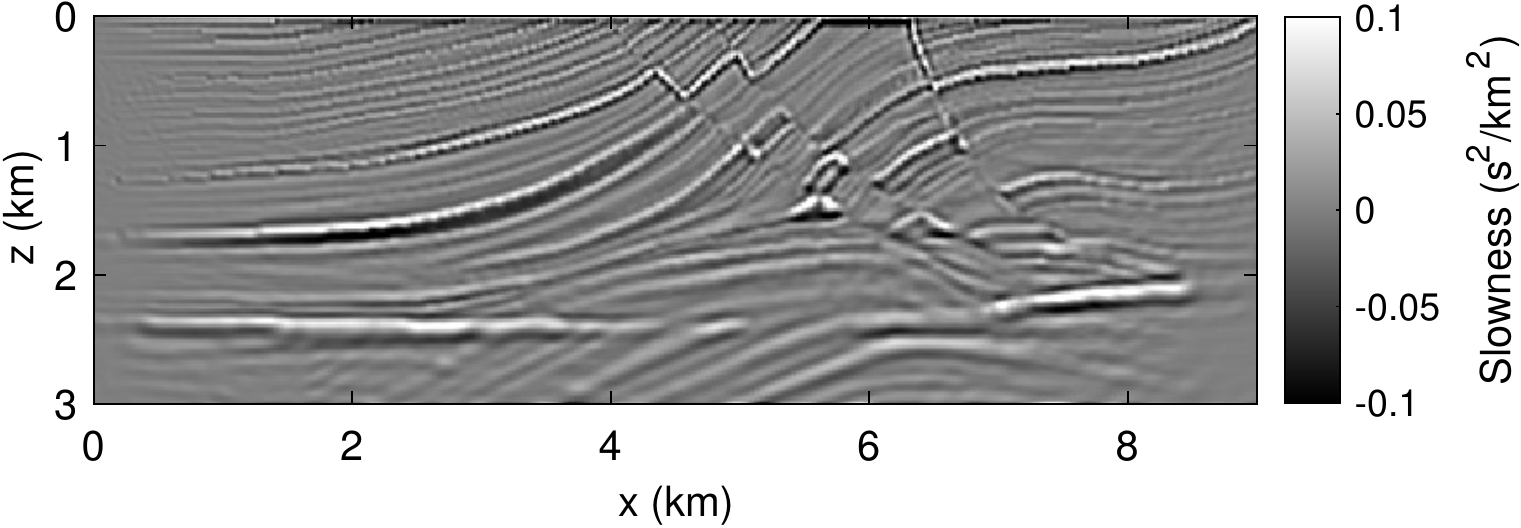}\label{fig:Mig-RTM-Lap}}
\caption{(a)~The smooth background velocity $m_0$ for RTM and LSRTM;~(b)~the true reflectivity; (c)~the original RTM result and (d)~the RTM image after Laplacian filtering.}~\label{fig:Mig-inv1}
\end{center}
\end{figure}

Our third example is to apply AA to LSRTM. We still use the Marmousi benchmark (\Cref{fig:Marm-true}) for illustration. The smooth background velocity is shown in~\Cref{fig:Mig-smooth}. We locate $80$ equally spaced wave sources (Ricker wavelet centered at 25 Hz) at $100$ m below the air-water interface. The entire workflow is similar to the FWI experiment except for a different forward problem and a different target. The size of the velocity model is $151$-by-$461$, and the spatial spacing is $20$ m. The total recording time is 4 seconds. The true reflectivity model is shown in~\Cref{fig:Mig-true}, and one iteration of RTM provides a crude subsurface image with unbalanced illumination, as seen in~\Cref{fig:Mig-RTM}. After the Laplacian filtering~\cite{zhang2009practical}, the migration artifacts are reduced, but the amplitude of the image is still incorrect, as seen by comparing the color bar of~\Cref{fig:Mig-RTM} and \Cref{fig:Mig-RTM-Lap} with the truth (\Cref{fig:Mig-true}). LSRTM aims to refine the image obtained by conventional RTM toward the true reflectivity. Therefore, we use the image obtained by the conventional RTM (\Cref{fig:Mig-RTM}) as the initial guess for inversion tests under the following four optimization methods: restarted GMRES, AA, L-BFGS, and the steepest descent. 


\begin{figure} 
\begin{center}
\subfloat[Restarted GMRES ($\cM=3$)]{ \includegraphics[width=0.7\textwidth]{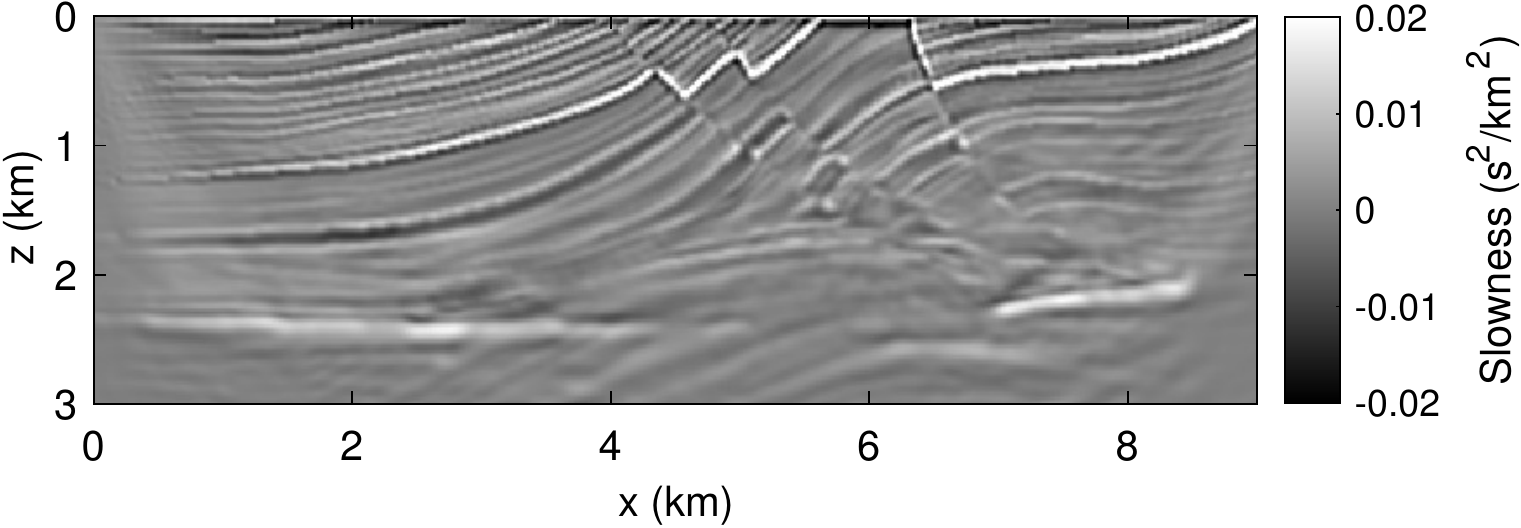}\label{fig:Mig-gmres}}\\
\subfloat[Anderson acceleration ($\cM=3$)]{ \includegraphics[width=0.7\textwidth]{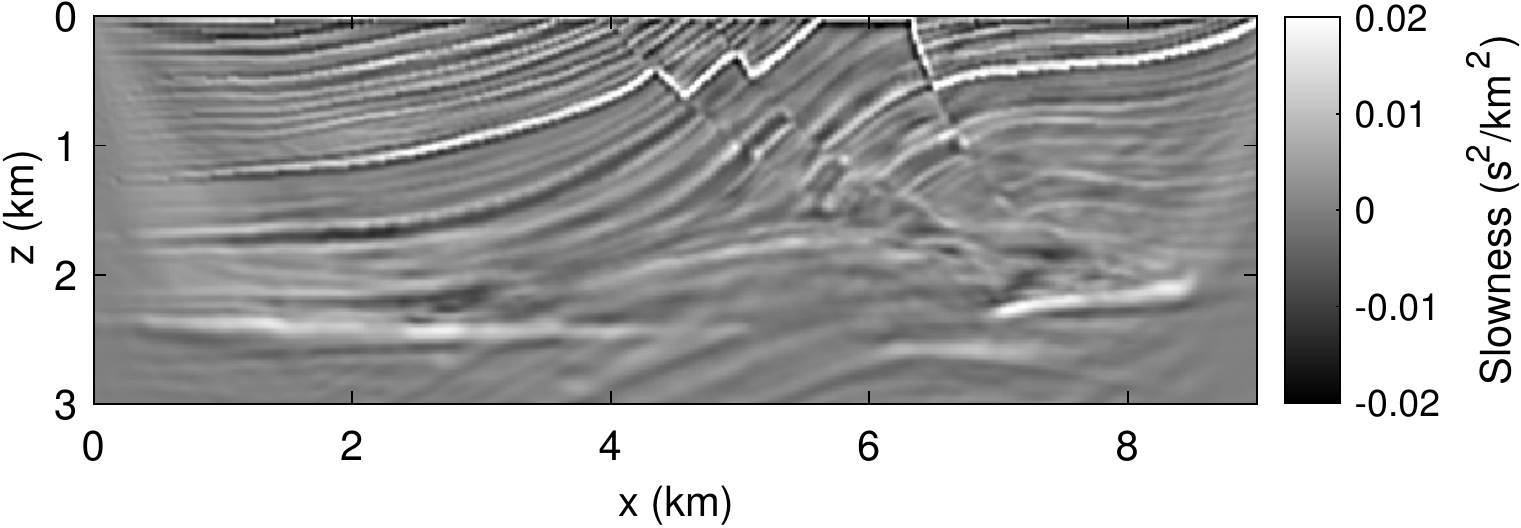}\label{fig:Mig-aa}}\\
\subfloat[L-BFGS ($\cM=3$)]{ \includegraphics[width=0.7\textwidth]{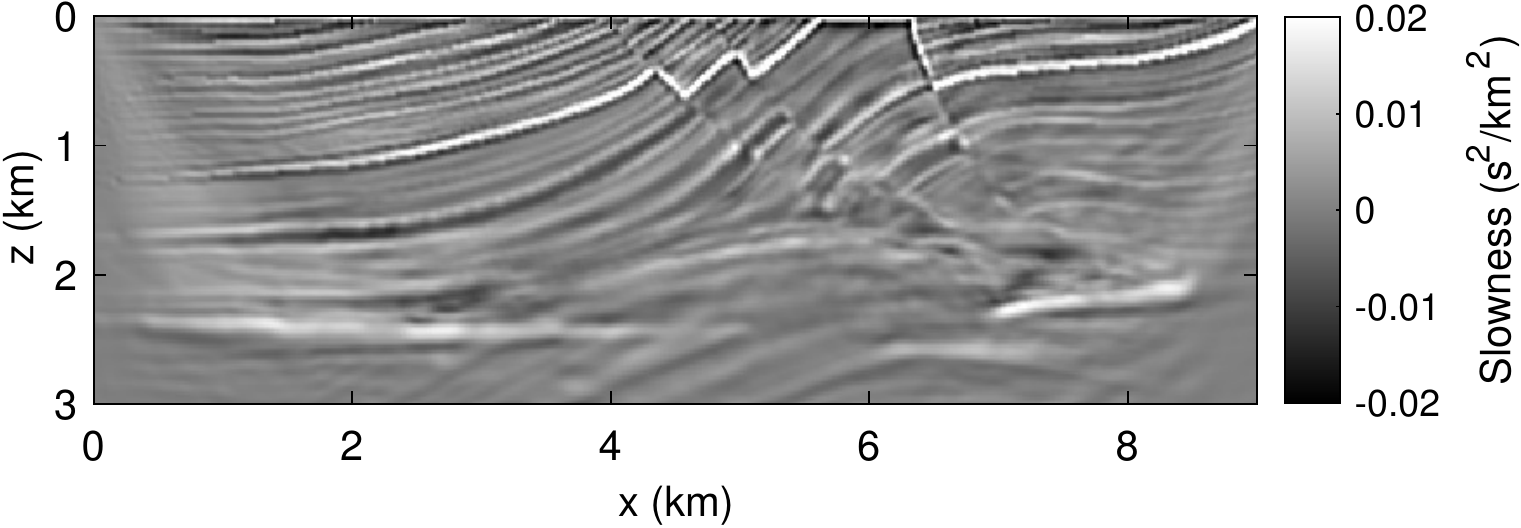}\label{fig:Mig-bfgs}}\\
\subfloat[Steepest descent]{ \includegraphics[width=0.7\textwidth]{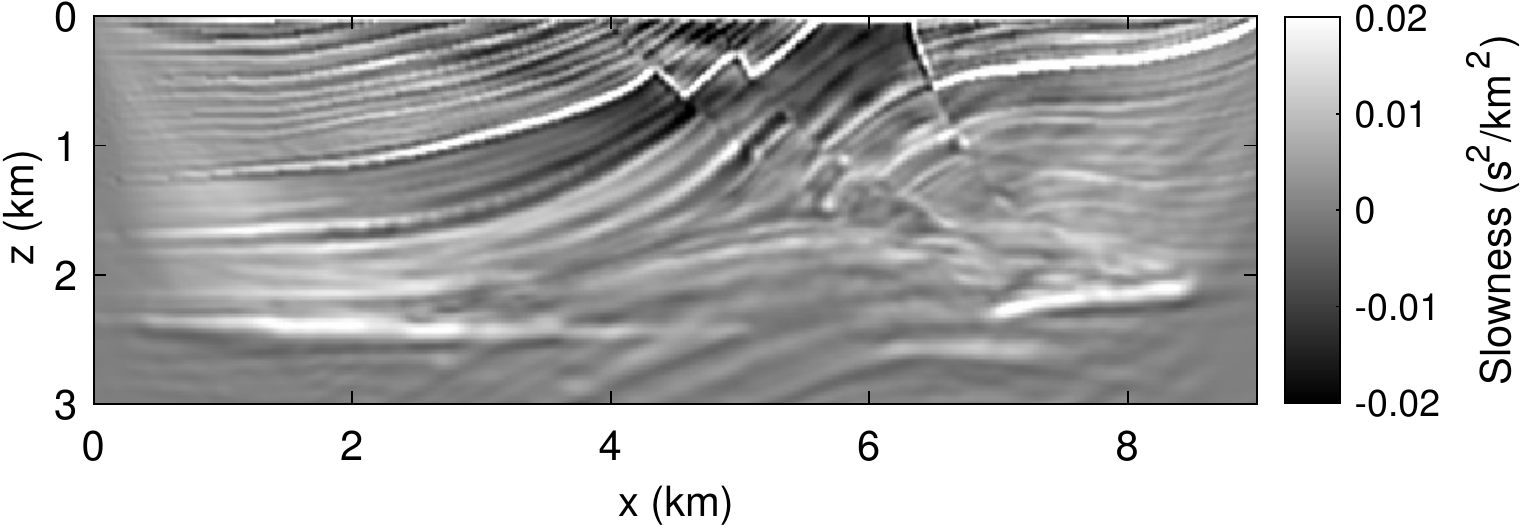}\label{fig:Mig-sd}}
\caption{LSRTM results by restarted GMRES, AA, L-BFGS, and the steepest descent method after $20$ iterations. The first three methods use the memory parameter $\cM=3$.}~\label{fig:Mig-inv2}
\end{center}
\end{figure}

Since GMRES is good at finding the solution for square linear systems, we reformulate the linear inverse problem that LSRTM aims to solve
\begin{equation}\label{eq:LSRTM1}
    Lm_r = d_r.
\end{equation}
We multiply both sides of~\Cref{eq:LSRTM1} by the Born operator $L$ and obtain
\begin{equation}\label{eq:LSRTM2}
    L^TLm_r = L^Td_r.
\end{equation}
Note the right-hand side of~\Cref{eq:LSRTM2} is nothing new but the migrated image after one step of RTM which we denote as $m_\text{RTM}$. Since $L^TL$ is a symmetric square matrix, we obtain a square linear system
\begin{equation}\label{eq:LSRTM3}
     A^L m_r = m_\text{RTM},
\end{equation}
where $A^L = L^TL$. Thus, we can use GMRES to find the solution of~\Cref{eq:LSRTM3}, which is also the solution of the original problem in~\Cref{eq:LSRTM1}. We use a restarted GMRES with memory parameter $\cM=3$. Therefore, at most three previous iterates are stored in memory when building up the Krylov space at each iteration. We further discuss the motivation and compare GMRES with AA in the next section. The final solution using the restarted GMRES is shown in~\Cref{fig:Mig-gmres} after 20 iterations.

The inverse problem or large-scale linear system that GMRES solves is an optimization-free formulation; see~\Cref{eq:LSRTM3}. Next, we return to the optimization formulation of the linear inverse problem~\eqref{eq:LSRTM} and use gradient-based methods to find the optimal $m^*_r$. Again, AA is applied to the steepest descent algorithm following~\Cref{alg:AA2}.~\Cref{fig:Mig-inv2} shows the inversion results after $20$ iterations. AA obtains an equally good image with the one by L-BFGS; see~\Cref{fig:Mig-aa} and~\Cref{fig:Mig-bfgs}. The memory parameter for both methods is $\cM=3$. It indicates that at most three iterates and their gradient vectors are stored in memory to compute the next iteration. AA also demonstrates a noticeable improvement in LSRTM compared with the steepest descent method, as shown in~\Cref{fig:Mig-sd}. All of the four methods improve the unbalanced illumination in the original RTM image (\Cref{fig:Mig-RTM}) and are closer to the true reflectivity (\Cref{fig:Mig-true}) compared to the filtered RTM image (\Cref{fig:Mig-RTM-Lap}).

\section{Discussions}
AA is a strategy proposed to speed up iterative schemes, particularly, fixed-point problems. In this paper, we modify the classical algorithm and apply it to accelerate iterative optimization algorithms, such as the gradient descent method. AA has intrinsic connections with GMRES and L-BFGS when solving specific types of problems. Therefore, we devote this section to detailed discussions and analysis regarding their differences and the potential benefits of AA. The ultimate goal is not to promote one method over another but to understand their roles in optimization problems better.

\subsection{Anderson Acceleration and GMRES}
The Generalized Minimal Residual Method (GMRES) is known as an iterative method to solve the square non-symmetric linear system~\cite{saad2003iterative}
\begin{equation}\label{eq:Axb}
Ax = b, \quad A\in \mathbb{C}^{n\times n}.
\end{equation}
We define the $k$-th Krylov subspace for this problem as
\begin{equation} \label{eq:Krylov}
K_k = K_k(A,r_0) = \text{span}\{r_0,Ar_0,\dots,A^{k-1}r_0\},
\end{equation}
where $r_0 = Ax_0 - b$ and $x_0$ is the initial guess. GMRES approximates the exact solution of $Ax=b$ by choosing the $k$-th iterate $x_k$ in the Krylov subspace $K_k$ such that
\begin{equation} \label{eq:GMRES}
x_k = \argmin_{x\in K_k} ||Ax_k - b||_2,
\end{equation}
which minimizes the Euclidean norm of the residual $r_k = Ax_k - b$. Since the current Krylov subspace is contained in the next subspace, 
$$ K_k =\text{span}\{x_0,x_1,\dots,x_{k-1}\} \subseteq K_{k+1} = \text{span}\{x_0,x_1,\dots,x_{k-1},x_{k}\},$$
the residual $r_k$ monotonically decreases as the number of iteration increases. There have been numerous variations and extensions of the method over the decades~\cite{saad2003iterative}.

In practice, it is not feasible to store all the previous iterates due to the limited machine memory. Instead, after $\cM$ iterations, one can treat the iterate $x_{\cM}$ as the new initial guess $x_0$, and constructs the new sequence of the Krylov subspace following~\eqref{eq:Krylov}. The method is well-known as the Restarted GMRES and often denoted as GMRES($\cM$). The Restarted GMRES may suffer from stagnation in convergence as the restarted subspace is often close to the earlier subspace for matrix $A$ with certain structures~\cite{embree2003tortoise}. Different from the Restarted GMRES, AA saves the memory by replacing the oldest iteration $x_0$ by the latest iterate $x_{\cM}$ if both methods share the same memory parameter $\cM$ and are applied to solve the square linear system~\eqref{eq:Axb}. The main differences between the two methods are illustrated in the following diagram.
\begin{eqnarray*}
\text{GMRES($\cM$):} &\quad \underbrace{\text{span}\{x_0,x_1,\dots,x_{\cM-1}\}}_{\text{the } \cM\text{-th subspace}} \longrightarrow \underbrace{\text{span}\{x_\cM\}}_{\text{the next subspace}}\\
\text{AA($\cM$):}& \quad \underbrace{\text{span}\{x_0,x_1,\dots,x_{\cM-1}\}}_{\text{the } \cM\text{-th subspace}} \longrightarrow \underbrace{\text{span}\{x_1,x_2,\dots,x_{\cM}\}}_{\text{the next subspace}}
\end{eqnarray*}

GMRES have been used in geophysical applications as a method of solving the \textit{forward} problems, which are also linear PDEs~\cite{erlangga2008iterative,calandra2012flexible}. One should note that the direct connections between AA and GMRES only apply when the fixed-point operator $G$ is based on the \textit{square linear problem}, as shown in~\eqref{eq:Axb}. So far, there has been no proved equivalence between AA applied to nonlinear operators and nonlinear GMRES. There have been numerical comparisons between AA of depth $\cM$ and GMRES($\cM$) in the literature~\cite{pratapa2016anderson,yang2020anderson}. Empirically, AA with memory parameter $\cM$ is observed to be more efficient than GMRES($\cM$) for certain nonlinear problems or linear problems with a non-positive definite matrix.

In FWI, we solve a nonlinear problem $F(x) = d$ where $F$ is the forward wave operator and $d$ is the observed data. Building Krylov spaces for such large-scale applications, which is theoretically an infinite-dimensional inverse problem, can be extremely costly. For LSRTM, we are solving a linear problem $Lm_r = d_r$ where $L$ is the Born operator, and $d_r$ is the observed scattering data. However, after discretization, $L$ becomes a tall skinny matrix. GMRES is not the best method for solving the linear system directly, and LSQR could be a better alternative. Therefore, we reformulate the problem of LSRTM so that it is suitable for GMRES. The results are presented in~\Cref{fig:Mig-inv2}. We remark that using AA for LSRTM accelerates the steepest descent algorithm, but GMRES for LSRTM, as it is done in this paper, is an optimization-free implementation.

\subsection{Anderson Acceleration and L-BFGS}
A standard optimization method used in geophysics is L-BFGS, where BFGS stands for the Broyden–Fletcher–Goldfarb–Shanno algorithm, and ``L'' indicates a variant of the algorithm with limited memory. We also have used this method for inversion tasks in the previous section. It belongs to the class of quasi-Newton methods, which is preferred when the full Newton’s Method is too time-consuming to apply. The Hessian matrix of a quasi-Newton method does not need to be computed explicitly at every iteration. Instead, an approximation $B_{k}$, which satisfies the following inverse secant condition, is used instead of the true inverse Hessian at the $k$-th iteration:
\begin{equation}~\label{eq:secant}
    B_{k} (J(p_{k}) - J(p_{k-1})) = p_{k} - p_{k-1}.
\end{equation}
Here, $p_{k-1}$ and $p_{k}$ are two consecutive iterates, and $J$ is the objective function that we aim to minimize. Quasi-Newton methods differ among each other in how to update $B_{k}$, the approximation to the inverse Hessian matrix. The BFGS algorithm follows two principles:~(1)~satisfy the inverse secant condition in~\Cref{eq:secant} and (2) be as close as possible to the approximation at the previous iteration. The latter is translated as
\begin{equation} \label{eq:bfgsFnorm}
    B_{k} = \argmin \limits_{B\in \mathbb{C}^{n\times n}} ||B - B_{k-1}||_F^2,
\end{equation}
where $\|\,\cdot\,\|_F$ denotes the matrix Frobenius norm, i.e., $\|A\|_F^2 = \sum_{i,j=1}^n|A_{ij}|^2$.
These conditions lead to explicit update schemes for BFGS and L-BFGS as a result of the famous Sherman--Morrison--Woodbury formula~\cite{nocedal2006numerical}. After replacing the inverse Hessian matrix in the Newton's method with the approximation $B_k$, the next iterate of the optimization problem is
\begin{equation} \label{eq:lbfgs_update}
    p_{k+1} = = p_k - B_kf_k = p_k + \eta B_{k}\G_k,
\end{equation}
where $\G_k$ denotes the gradient of the objective function $J(p)$ at $p=p_k$.

While the BFGS algorithm is a type of secant method, it is worth addressing that AA is equivalent to a \textit{multi-secant} method~\cite{fang2009two}. If one sticks with $\ell^2$-based AA regarding~\Cref{eq:AAopt} as we do in this paper, the update formula~\eqref{eq:AAupdate0} (with $\beta_k=1$) can be re-written as: 
\begin{equation} \label{eq:aa_secant_update}
    p_{k+1} = p_{k}  - S_k f_k = p_{k}  + \eta S_k \G_k,
\end{equation}
where $S_k \in \mathbb{C}^{n\times n}$ is the solution to the following constrained optimization problem:
\begin{equation} \label{eq:aaFnorm}
S_k = \argmin \limits_{S \in\mathbb{C}^{n\times n}} \|S +I\|^2_{F}
\end{equation} 
subject to the multi-secant condition
\begin{equation}\label{eq:multisecant}
    S_k D_k = P_k.
\end{equation}
Denoting $\Delta p_i = p_{i+1}-p_i$, $\Delta f_i = f_{i+1}-f_i$, matrices $P_k$ and $D_k$ are defined as 
\begin{equation} 
P_k = [\Delta p_{k-m},\ldots,\Delta p_{k-1}]\in\mathbb{C}^{n\times m},\quad D_k = [\Delta f_{k-m_k},\ldots,\Delta f_{k-1}]\in\mathbb{C}^{n\times m}.
\label{eq:XkDk}
\end{equation}
We have used the fact that $f_k=G(p_k)-p_k =- \eta\,\G_k$ in~\Cref{eq:aa_secant_update}.

Once we have rewritten the AA algorithm, it is not hard to recognize the similarities between the update formula of AA and the one of the L-BFGS method; see~\Cref{eq:lbfgs_update} and~\Cref{eq:aa_secant_update}. The key difference is that matrices $S_k$ and $B_k$ are constructed under different principles, although both can be regarded as approximations to the inverse Hessian matrix in the full Newton's method. For L-BFGS, $B_k$ satisfies the secant condition~\eqref{eq:secant} while minimizing~\eqref{eq:bfgsFnorm}. For AA, $S_k$ satisfies the \textit{multi-secant} condition~\eqref{eq:multisecant} while minimizing~\eqref{eq:aaFnorm}. Both methods are faster than the steepest descent algorithm and have been proved to have superlinear convergence.

The original BFGS algorithm stores a dense $n$-by-$n$ approximation to the inverse Hessian matrix, where $n$ is the number of variables. Besides, each BFGS iteration has a cost of $\bO(n^2)$ arithmetic operations. The idea of L-BFGS is to restrict the use of all iterations in the history to the latest $\cM$ iterates; $\cM$ is a parameter of the L-BFGS algorithm that can be chosen a priori. Since the earlier iterates often carry little information about the curvature of the current iterate, the change from BFGS to L-BFGS is expected to have minimal impacts on the convergence rate. The L-BFGS method has a linear memory requirement in terms of the number of variables. It is particularly well-suited for large-scale optimization problems such as FWI and LSRTM.

Similar to AA, L-BFGS maintains a history of the past $\cM$ iterates and their gradients. Again, $\cM$ is often chosen to be small. The connections are illustrated by the following diagram where $\G_k$ denotes the gradient of the objective function at $p = p_k$.
\begin{eqnarray*}
\text{L-BFGS($\cM$):} &\quad \underbrace{\{p_0,\dots,p_{\cM-1}, \G_0,\dots,\G_{\cM-1}\}}_{\text{construct } B_{\cM-1}} \longrightarrow \underbrace{\{p_1,\dots,p_{\cM}, \G_1,\dots,\G_\cM\}}_{\text{construct } B_{\cM}}\\
\text{AA($\cM$):} &\quad
\underbrace{\{p_0,\dots,p_{\cM-1},\G_0,\dots,\G_{\cM-1}\}}_{\text{construct } S_{\cM-1}} \longrightarrow \underbrace{\{p_1,\dots,p_{\cM},\G_1,\dots,\G_\cM\}}_{\text{construct } S_{\cM}}
\end{eqnarray*} 
We remark that the relationship above is valid only when AA is used to accelerate the steepest descent algorithm as what we propose in this paper. The most expensive part of the L-BFGS algorithm is the inverse Hessian update, and the most costly step of AA is to compute the optimal coefficients. In terms of computational cost, the low-rank QR update for AA, the inverse Hessian update for L-BFGS, and the Restarted GMRES all take $\bO(n\cM)$ floating-point operations (flops), if implemented optimally. The minimum memory requirements for all three methods are also $\bO(n\cM)$, considering the memory parameter for all the methods is $\cM$, and the size of the unknown is $n$.


\subsection{Performance comparison}
We have seen the theoretical connections among AA, the Restarted GMRES, and the L-BFGS algorithm in the previous two subsections. Although under the same memory parameter $\cM$, all the methods have the same order of computational cost and memory demand, AA could be advantageous in the following two aspects.

First, when the iteration number is bigger than the memory parameter $\cM$, AA always uses the last $\cM$ iterates to construct the solution for the next iteration. On the other hand, the Restarted GMRES throws away all the $\cM$ iterates and restarts from zero once the restart window is reached. Thus, AA can use more information in the optimization than the Restarted GMRES in most iterations. The advantage is reflected in our numerical examples for LSRTM; see~\Cref{fig:Mig-gmres} and~\Cref{fig:Mig-aa}. 

One can only compare AA and GMRES when both methods are used to find the solution $x$ for problems in the form of $Ax = b$, where $ A $ is a square matrix. The flexibility of GMRES fades for rectangle matrices and highly nonlinear problems, while AA can be applied to all types of linear and nonlinear problems as long as they are written as fixed-point operators. Some study has shown that the fixed-point operator does not need to be a contraction for AA to converge, although it is necessary for Picard iteration~\cite{pollock2019anderson}.

Second, unlike the Restarted GMRES, the L-BFGS algorithm uses all the last $\cM$ iterates to compute the next iteration. Although both AA and L-BFGS exploit all the information available in storage, the two methods approximate the inverse Hessian matrix in slightly different ways: the former is a \textit{multi-secant} method while the latter is a secant method. We observe in \Cref{fig:inv-aa} and~\Cref{fig:inv-bfgs} that FWI using AA spends less time searching for an appropriate step size than the inversion using the L-BFGS method on average for each iteration. It helps AA achieve better inversion results than L-BFGS under the same number of total gradient evaluations. The inverse Hessian matrix approximated by AA satisfies the secant condition not only for the latest iteration but also for the previous $\cM$ iterations. This implicitly enforces the connections among the neighbor iterates to avoid unstable descent directions far from the curvature of the basin of attraction.

Although it is expected that the bigger the memory parameter $\mathcal{M}$ for AA, the better the performance, empirical studies have shown that a relatively small $\mathcal M$, commonly ranging from $3$ to $20$, is often good enough to speed up the convergence of the fixed-point iteration without a significant toll on the machine memory and the computation cost. In the experiments, the choice of the memory parameter for AA could follow similar principles as choosing the memory depth for L-BFGS and the restarted window for GMRES.


\section{Conclusion}

Anderson acceleration for seismic inversion treats the method of steepest descent as a fixed-point operator. It speeds up the convergence by linearly combining a list of the previous iterates in an optimized way. The computational cost of implementing AA mainly comes from the 1D optimization for the weights. It is thus easy to add AA to existing optimization algorithms. As shown in the paper, AA outperforms the steepest descent method and can also be considered an alternative to L-BFGS. AA is equivalent to a \textit{multi-secant} method, while the L-BFGS algorithm is derived by the secant method. Being computationally attractive, AA is a promising optimization algorithm for FWI and LSRTM.

\section{Acknowledgements}
The author thanks Dr.~Alan Richardson and other anonymous reviewers for constructive suggestions for the paper. This work is supported in part by the National Science Foundation through grant DMS-1913129. 

\bibliographystyle{cpam}
\bibliography{AA4FWI}

\end{document}